\documentclass[11pt]{article}
\usepackage{mathrsfs}
\usepackage{stmaryrd}
\usepackage{amsfonts}
\usepackage{amsfonts,amsmath,amssymb,amscd}
\usepackage{shadow}
\usepackage{graphicx}
\usepackage{color}
\usepackage{longtable}
\allowdisplaybreaks
\newtheorem{theo}{Theorem}[section]

\newtheorem{prop}[theo]{Proposition}
\newtheorem{coro}[theo]{Corollary}

\def\qed{\hfill \rule{4pt}{7pt}}

\def\pf{\noindent {\it{Proof.} \hskip 2pt}}

\def \inv{\mathop{\rm inv }}
\def \maj{\mathop{\rm maj }}

\def \des{\mathop{\rm des }}
\def \Des{\mathop{\rm Des }}

\def \Z{\mathop{\rm Z }}

\def \C{\mathscr{C}}

\begin{document}
\begin{center}
{\LARGE\bf  Han's Bijection via Permutation Codes}
\end{center}

\begin{center}
William Y.C. Chen$^1$, Neil J.Y. Fan$^2$, Teresa X.S. Li$^3$

Center for Combinatorics, LPMC-TJKLC\\
Nankai University, Tianjin 300071, P.R. China

$^1$chen@nankai.edu.cn, $^2$fjy@cfc.nankai.edu.cn,
$^3$lxs@cfc.nankai.edu.cn.

\end{center}

\begin{abstract}
We show that Han's bijection when restricted to permutations can be
carried out in terms of the cyclic major code and the cyclic inversion
code. In other words,
 it maps a permutation $\pi$ with a cyclic major code $(s_1, s_2, \ldots, s_n)$ to a
permutation $\sigma$ with a cyclic inversion code $(s_1,s_2, \ldots,
s_n)$. We also show that the fixed points of Han's map can be
characterized by the strong fixed points of Foata's second
fundamental transformation. The notion
 of strong fixed points is related to  partial
 Foata maps introduced by Bj\"{o}rner and Wachs.

\end{abstract}



\section{Introduction}

In his combinatorial proof of the fact that the $\Z$-statistic
introduced by Zeilberger and Bressoud \cite{Bress-Zeil} is Mahonian,
Han \cite{Han3} constructed a Foata-style bijection on words which
maps the major index onto the $\Z$-statistic. Since the
$\Z$-statistic and the inversion number coincide when restricted to
permutations,  Han's bijection maps  the major
index to the inversion number for permutations.  Let $H$ denote
Han's bijection when restricted to permutations. Throughout this
paper, by Han's bijection we always mean
 the map $H$. We shall show that the map
$H$ can be carried out by the cyclic major code and the cyclic
inversion code.

The cyclic major code of a
 permutation can be described in terms of cyclic intervals,
a notion also introduced by Han \cite{Han2} in his study of the
joint distribution of the excedance number and Denert's statistic.
It should be noted that the cyclic inversion code in the context of
this paper is the classical Lehmer code, but  the cyclic major code
is different from the well-studied major code
as introduced by Rawlings \cite{Rawlings}, see also Foata and Han \cite{Foata-Han},
Skandera \cite{Skandera}, Dzhumadil'daev \cite{Dzhumadil},  and  Han \cite{Han4}.

 Using the code representation,
we show that the fixed points of Han's map can be characterized by the strong fixed points of Foata's second fundamental transformation.
The notion of strong fixed points is related to partial Foata maps introduced by Bj\"{o}rner and Wachs \cite{bjoner}.

Let us give an overview of the background and definitions. Let $X=\{1^{m_1},2^{m_2},\cdots,k^{m_k}\}$ be a multiset with $m_i$
$i$'s and $m_1+m_2+\cdots+m_k=n$. The set of    rearrangements of
 $X$ is denoted by $R(X)$.
When $m_1=m_2=\cdots=m_k=1$,
$R(X)$ reduces to the set $S_n$ of permutations on $[n]$.
For a word
$w=w_1w_2\cdots w_n\in R(X)$, the  descent set  $\Des(w)$, the descent
number $\des(w)$, the  major index  $\maj(w)$,  the inversion number $\inv(w)$
and the  $\Z$-statistic  $\Z(w)$ are defined by
\begin{align*}
\Des(w)&=\{i|1\leq i\leq n-1,w_i>w_{i+1}\},\\[5pt]
\des(w)&=\#\Des(w),\\[5pt]
\maj(w)&=\sum_{i\in \Des(w)}i,\\[5pt]
\inv(w)&=\#\{(i,j)|1\leq i<j\leq n, w_i>w_j\},\\[5pt]
\Z(w)&=\sum_{i<j} \maj(w_{ij}),
\end{align*}
where $w_{ij}$ is a word obtained from $w$ by deleting all elements
except $i$ and $j$. For example, let  $w=211324314\in
R(1^3,2^2,3^2,4^2)$. We have $\Des(w)=\{1,4,6,7\}, \des(w)=4,
\maj(w)=18, \inv(w)=9$, and $\Z(w)$ can be computed as follows
\[ \maj(21121)+\maj(11331)+\maj(11414)+\maj(2323)+
\maj(2244)+\maj(3434)=16.\]

A statistic is said to be
Mahonian on $R(X)$ if it has the same distribution as the major index on
$R(X)$.  MacMahon  \cite{Macmahon1, Macmahon2}  introduced the major index   and
 proved that the major index is
equidistributed with the inversion number for $R(X)$.
Foata \cite{Foata} found a combinatorial
proof of this classical fact  by constructing a bijection $\Phi$,
called the second fundamental transformation, which maps the major index
 to the inversion number, namely,
\begin{align*}
\maj(w)=\inv(\Phi(w)) \ \ {\rm for}\  {\rm any}\  w\in R(X).
\end{align*}

For completeness, we give a brief description of Foata's bijection \cite{Foata},
 see also Haglund \cite{Haglund}, Lothaire \cite{Lothaire}.
Let $w=w_1w_2\cdots w_n$ be a word on a multiset $X$ as defined
above, and let $x$ be an element in $X$.  If $w_n\leq x$, the
$x$-factorization of $w$ is defined as $w=v_1b_1\cdots v_pb_p$,
where each $b_i$ is less than or equal to $x$, and every element in
$v_i$
 is greater than $x$. Note that $v_i$ is allowed to be empty.
Similarly, when $w_n>x$, the $x$-factorization of $w$ is
defined as $w=v_1b_1\cdots v_pb_p$, where each $b_i$ is  greater than
$x$, and every element in $v_i$ is  less than or equal to $x$. In either case, set
\[ \gamma_x(w)=b_1v_1\cdots b_pv_p, \quad
w'=w_1w_2\cdots w_{n-1}.\] Then the second fundamental
transformation $\Phi$ can be defined recursively by setting $\Phi(a)=a$ for each
$a\in X$ and setting
\begin{align*}
\Phi(w)=\gamma_{w_n}(\Phi(w'))\cdot w_n
\end{align*}
if $w$ contains more than one element.

As an extension of the theorem of MacMohan, Bj\"orner and Wachs \cite{bjoner} considered
the problem of finding subsets $U$ of $S_n$ for which the major index and
inversion number are equidistributed.
They   introduced the $k$-th partial Foata
bijection $\phi_k\colon S_n\longrightarrow S_n $ for $1\leq k \leq n$. Let $\sigma=\sigma_1\sigma_2\cdots\sigma_n\in S_n$. Define $
\phi_1(\sigma)=\sigma$ and for $k>1$ define
\[ \phi_k(\sigma)=\gamma_{\sigma_k}(\sigma_1\sigma_2\cdots\sigma_{k-1})
\cdot\sigma_k\sigma_{k+1}\cdots\sigma_n.
\]
It is easily seen that
 \[ \Phi=\phi_n\circ\phi_{n-1}\cdots\circ\phi_1.\]
 A subset $U$ of $S_n$ is said to be a
strong Foata class if \[ \phi_{k}(U)=U\] for $1\leq k \leq n$. A permutation $\sigma$ is said to be  a strong  fixed point of Foata's map if
 \[ \phi_{k}(\sigma)=\sigma\] for  $1\leq k\leq  n.$  As will be seen,
the strong fixed points of Foata's map is closely related to
the fixed points of Han's map.

The paper is organized as follows. In Section 2, we recall   the
construction of Han's map, and give a description of the cyclic
major code and the cyclic inversion code. Then we give a reformulation of Han's map in terms of
these two codes.  In Section 3,  we give a characterization of the fixed points of Han's map
$H$. It turns out that a permutation is fixed by $H$ if and only if it is a
strong fixed points of Foata's map $\Phi$.

\section{Han's bijection via permutation codes}

In this section, we are  concerned with  a  reformulation  of Han's
bijection for permutations in terms of the cyclic major code and the
cyclic inversion code.
 For completeness, let us give an overview of
the  map $H$.

Let $x\in[n]$ and
$\sigma=\sigma_1\sigma_2\cdots\sigma_{n-1}$ be a permutation on
$\{1,2,\cdots,x-1,x+1,\cdots,n\}.$ Define
$C^x(\sigma)$ as $\tau_1\tau_2\cdots\tau_{n-1}$, where $\tau_i=\sigma_i-x({\rm mod}\ n)$, i.e.,
\begin{align*}
\tau_i=\left\{
         \begin{array}{ll}
           \sigma_i-x+n, & \hbox{if  $\sigma_i<x$}; \\[5pt]
           \sigma_i-x, & \hbox{if $\sigma_i>x$},
         \end{array}
       \right.
\end{align*}
and define $C_x(w)$  as the
standardization of $\sigma$, i.e.,
$C_x(w)=\nu_1\nu_2\cdots\nu_{n-1}\in S_{n-1}$ with
\begin{align*}
\nu_i=\left\{
        \begin{array}{ll}
          \sigma_i, & \hbox{if $\sigma_i<x$;} \\[5pt]
          \sigma_i-1, & \hbox{if $\sigma_i>x$.}
        \end{array}
      \right.
\end{align*}
Evidently, both $C^x$ and $C_x$ are bijections between permutations on
$\{1,2,\cdots,x-1,x+1,\cdots,n\}$ and $S_{n-1}$.
So $(C^x)^{-1}$ and $(C_x)^{-1}$ are well defined.
Han's bijection $H$  can be defined by $H(1)=1$ and
\begin{align*}
H(\sigma)=C_{\sigma_n}^{-1}(H(C^{\sigma_n}(\sigma')))\cdot\sigma_n,
\end{align*}
where $\sigma\in S_n$ with $n>1$ and $\sigma^{'}=\sigma_1\sigma_2\cdots\sigma_{n-1}$.

We proceed to give the definition of cyclic intervals. Let $X=\{1^{m_1},2^{m_2},\ldots,k^{m_k}\}$ be a multiset. For  $x,y \in X$, the
cyclic interval $\rrbracket x,y\rrbracket$ is defined by Han
\cite{Han2} as
\begin{align*}
\rrbracket x,y\rrbracket=\left\{
                           \begin{array}{ll}
                             \{z|z\in [k],x<z\leq y\}, & \hbox{if $x\leq y$;} \\[5pt]
                             \{z|z\in [k], z>x\ \ or\ \ z\leq y\}, & \hbox{otherwise.}
                           \end{array}
                         \right.
\end{align*}
Set $\rrbracket x, \infty\rrbracket=\{z|z\in [k],z>x\}.$

For any word $w=w_1w_2\cdots w_n$  on  $X$ and $1\leq i\leq n$, define
\begin{align*}
t_i(w)&=\#\{j|1\leq j\leq i-1, w_j\in \rrbracket w_i,
\infty\rrbracket\},
\end{align*} and
\begin{align*}
s_i(w)&=\#\{j|1\leq j\leq i-1, w_j\in
\rrbracket w_i, w_{i+1}\rrbracket\},
\end{align*}
where $w_{n+1}=\infty,$ and $\# S$ stands for the cardinality of a set $S$.

For example, let $w=312432143$. Then
\[(t_1(w),t_2(w),\cdots,t_9(w))=(0, 1, 1, 0, 1, 3, 5, 0, 2),\]
and
\[(s_1(w),s_2(w),\cdots,s_9(w))=(0, 0, 1, 3, 3, 4, 5, 6, 2).\]
The notion of cyclic intervals plays an important role in
 the proof of the fact that the bi-statistic
$(\rm
{exc
}, \rm{Den})$ is equidistributed with  $(\des, \maj)$ on $R(X)$, where  ${\rm exc}$ is the excedance number
and  ${\rm Den}$ is the  Denert's statistic, see Denert \cite{Denert},  Foata and Zeilberger \cite{Foata-Zeil}, and Han \cite{Han2}.

We proceed to give the definition of the cyclic  major code also in
terms of cyclic intervals. Meanwhile, the traditional inversion code
can  be described in this way. Let \[ E_n=\{(a_1,a_2,\cdots,a_n)\in
Z^n|0\leq a_i\leq i-1,i=1,2,\cdots,n\}.\] Keep in mind that the
above definitions of   $t_i(\sigma)$ and $s_i(\sigma)$  apply to
permutations. It is well known that the map $I\colon
S_n\longrightarrow E_n$ defined by \[
\sigma\longmapsto(t_1(\sigma),t_2(\sigma),\cdots,t_n(\sigma))
\]
is a bijection known as the Lehmer code, which is often referred to as the inversion code.
Note that
\[ \sum_{i=1}^nt_i(\sigma)=\inv(\sigma).\] On the other hand, it is
easy to see that the  map
$M\colon S_n\longrightarrow E_n$ defined by
\[ \sigma\longmapsto(s_1(\sigma),s_2(\sigma),\cdots,s_n(\sigma))
\]
is also a bijection. We call $M(\sigma)$  the cyclic major code of
$\sigma$.  To recover $\sigma$ from its cyclic major code
$(s_1,s_2,\ldots,s_n)$, first let $\sigma_n=n-s_n$. Suppose that
$\sigma_{k+1},\ldots,\sigma_n$ have been determined by
$s_{k+1},\ldots,s_n$. Then delete the elements in the sequence
\[
 \sigma_{k+1},\sigma_{k+1}-1,\ldots,1,n,(n-1),\ldots,(\sigma_{k+1})+1
\]
that are equal to $\sigma_j$ for some $j\geq k+1$ and set $\sigma_k$ to be the $(s_{k}+1)$-th element in the resulting sequence.   It has been shown by Han \cite{Han2} that
\[ \sum_{i=1}^ns_i(\sigma)=\maj(\sigma).\]
For example, $I(38516427)=(0,0,1,3,1,3,5,1)$ and $M(38516427)=(0,1,1,2,3,4,4,1)$.

The relation between these two codes is described below.

\begin{prop}\label{tranfor} Let $\sigma\in S_n$. Suppose that
$I(\sigma)=(t_1,t_2,\cdots,t_n)$ and
$M(\sigma)=(s_1,s_2,\cdots,s_n)$. Then we have $s_n=t_n$, and for $1\leq i<
n$, $s_i=t_i-t_{i+1}({\rm mod}\ i)$, that is,
\[ s_i=\left\{
         \begin{array}{ll}
           t_i-t_{i+1}, & \hbox{if $t_i\geq t_{i+1}$;} \\[5pt]
           t_i-t_{i+1}+i, & \hbox{if $t_i< t_{i+1}$.}
         \end{array}
       \right.\]
\end{prop}

\pf It is clear that $s_n=t_n$. For $1\leq i\leq n-1$, by the definition of
$t_i(\sigma)$, we see that $t_i\geq t_{i+1}$ if and only if
$\sigma_i<\sigma_{i+1}$. In this case,
\begin{align*}
s_i&=\#\{j| 1\leq j\leq i-1,\sigma_i<\sigma_j<\sigma_{i+1}\}\\[5pt]
&=\#\{j| 1\leq j\leq i-1,\sigma_i<\sigma_j\}-\#\{j| 1\leq j\leq
i-1,\sigma_j>\sigma_{i+1}\}\\[5pt]
&=t_i-t_{i+1}.
\end{align*}

If $t_i<t_{i+1}$, then $\sigma_i>\sigma_{i+1}$. Hence
\begin{align*}
s_i&=\#\{j| 1\leq j\leq i-1,\sigma_i<\sigma_j\ {\rm or }\
\sigma_j<\sigma_{i+1}\}\\[5pt]
&=\#\{j| 1\leq j\leq i-1,\sigma_i<\sigma_j\}+\#\{j| 1\leq j\leq
i-1,\sigma_j<\sigma_{i+1}\}\\[5pt]
&=t_i+\#\{j| 1\leq j\leq i-1,\sigma_j<\sigma_{i+1}\}\\[5pt]
&=t_i+i-\#\{j| 1\leq j\leq i,\sigma_j>\sigma_{i+1}\}\\[5pt]
&=t_i-t_{i+1}+i.
\end{align*}
This completes the proof. \qed

The following theorem states that  Han's bijection $H$ can be
carried out in terms of the cyclic major code and the   inversion code.

\begin{theo}\label{th1}
For each $n\geq 1$,  we have
\[H=I^{-1}\circ M .\]
In other words, $H$ is a bijection on $S_n$ with the property that
\[M(\sigma)=I(H(\sigma)).\]
\end{theo}
\pf  We use
induction on $n$. For $n=1$, the theorem is obvious. Assume that $n>1$. Let $\sigma=\sigma_1\sigma_2\cdots\sigma_n\in S_n$, and let $M(\sigma)=(s_1(\sigma),s_2(\sigma),\cdots,s_n(\sigma))$. By
definition, $s_n(\sigma)=\#\{\sigma_n+1,\sigma_n+2,\cdots,n\}=n-\sigma_n$.
By the construction of $H$, we have
\[H(\sigma)=C_{\sigma_n}^{-1}[H(C^{\sigma_n}(\sigma')]\cdot\sigma_n,\]
which implies
$t_n(H(\sigma))=n-\sigma_n$.
Since the standardization of a permutation  preserves  the relative order, we find that
\[I(C_{\sigma_n}(\sigma_1\sigma_2\cdots\sigma_{n-1}))=(t_1(\sigma),t_2(\sigma),\cdots,t_{n-1}(\sigma)).\]
By induction, it suffices to show that
\begin{align}\label {eq1}
M(C^{\sigma_n}(\sigma_1\sigma_2\cdots\sigma_{n-1}))=(s_1(\sigma),s_2(\sigma),\cdots,s_{n-1}(\sigma)).
\end{align}
Suppose that
$C^{\sigma_n}(\sigma_1\sigma_2\cdots\sigma_{n-1})=\tau_1\tau_2\cdots\tau_{n-1}$.
For the sake of presentation, let  $\tau_n=\infty$. For $1\leq i \leq n-1$ and $1\leq k\leq i-1$,  we claim that
$\sigma_k\in \rrbracket \sigma_i,\sigma_{i+1}\rrbracket$ if and only
if $\tau_k\in \rrbracket \tau_i,\tau_{i+1}\rrbracket$. If it is true,
then  \eqref{eq1} follows immediately. This claim can be verified as
follows.

(1) If $i\neq n-1$, there are two cases each of which
has three subcases, namely,
\begin{align*}
&(1a)\ \ \sigma_n<\sigma_i<\sigma_{i+1};\\[5pt]
&(1b)\ \ \sigma_i<\sigma_n<\sigma_{i+1};\\[5pt]
&(1c)\ \ \sigma_i<\sigma_{i+1}<\sigma_n;
\end{align*}
and
\begin{align*}
&(2a)\ \ \sigma_n>\sigma_{i}>\sigma_{i+1};\\[5pt]
&(2b)\ \ \sigma_i>\sigma_{n}>\sigma_{i+1};\\[5pt]
&(2c)\ \ \sigma_i>\sigma_{i+1}>\sigma_{n}.
\end{align*}
We only give the proof of case (1b), the other cases can be justified
by the same argument. Let us assume that $\sigma_i<\sigma_n<\sigma_{i+1}$. By
definition, $\tau_i=n+\sigma_i-\sigma_n,
\tau_{i+1}=\sigma_{i+1}-\sigma_{n}$, so we have
$\tau_{i+1}<\tau_i$.
Suppose that $\sigma_k\in\rrbracket \sigma_i,\sigma_{i+1}\rrbracket$.
Then we deduce that $\sigma_i<\sigma_k<\sigma_{i+1}$ and
\[\tau_k=\left\{
           \begin{array}{ll}
             \sigma_k-\sigma_n+n, & \hbox{if $\sigma_k<\sigma_n<\sigma_{i+1}$;} \\[5pt]
             \sigma_k-\sigma_n, & \hbox{if $\sigma_i<\sigma_n<\sigma_{k}$.}
           \end{array}
         \right.
\]
If $\sigma_k<\sigma_n<\sigma_{i+1}$, then
$\tau_k=\sigma_k-\sigma_n+n>\sigma_i-\sigma_n+n=\tau_i$, it follows that
$\tau_k\in\rrbracket \tau_i,\tau_{i+1}\rrbracket$; if
$\sigma_i<\sigma_n<\sigma_{k}$, then
$\tau_k=\sigma_k-\sigma_n<\sigma_{i+1}-\sigma_n=\tau_{i+1}$, which implies $\tau_k\in\rrbracket \tau_i,\tau_{i+1}\rrbracket$. Conversely,
if  $\tau_k\in\rrbracket \tau_i,\tau_{i+1}\rrbracket$, then we deduce that
$\tau_k>\tau_i$ or $\tau_k<\tau_{i+1}$. Assume that $\sigma_k\notin \rrbracket
\sigma_i,\sigma_{i+1}\rrbracket$, then we have $\sigma_k<\sigma_i$ or
$\sigma_k>\sigma_{i+1}$. Consequently,
\[\tau_k=\left\{
           \begin{array}{ll}
             \sigma_k-\sigma_n+n, & \hbox{if $\sigma_k<\sigma_i<\sigma_n$;} \\[5pt]
             \sigma_k-\sigma_n, & \hbox{if $\sigma_k>\sigma_{i+1}>\sigma_n$.}
           \end{array}
         \right.
\]
If $\sigma_k<\sigma_i$, then
$\tau_k=\sigma_k-\sigma_n+n<\sigma_i-\sigma_n+n=\tau_i$. However,
 $\tau_k=\sigma_{k}+n-\sigma_{n}>\sigma_{i+1}-\sigma_n=\tau_{i+1}$, which is a
contradiction. If $\sigma_k>\sigma_{i+1}$, then
$\tau_k=\sigma_k-\sigma_n>\sigma_{i+1}-\sigma_n=\tau_{i+1}$, but now
$\tau_k=\sigma_k-\sigma_n<\sigma_{i}-\sigma_{n}+n=\tau_i$, a
contradiction too. So we reach the conclusion that  $\sigma_k\in \rrbracket
\sigma_i,\sigma_{i+1}\rrbracket$.

(2) If $i=n-1$, there are two cases, namely $\sigma_{n-1}>\sigma_n$
and  $\sigma_{n-1}<\sigma_n.$
For the first case, by  definition we have $\tau_{n-1}=\sigma_{n-1}-\sigma_n$. It follows that
\begin{align*}
\sigma_k\in \rrbracket \sigma_{n-1},\sigma_{n}\rrbracket&\Rightarrow
\sigma_k>\sigma_{n-1}\ \ {\rm or}\ \ \sigma_k<\sigma_n\\[5pt]
&\Rightarrow \tau_k=\left\{
          \begin{array}{ll}
            \sigma_k-\sigma_{n}, & \hbox{if $\sigma_k>\sigma_{n-1}$;} \\[5pt]
            \sigma_k-\sigma_{n}+n, & \hbox{if $\sigma_k<\sigma_{n}$.}
          \end{array}
        \right.\\[5pt]
&\Rightarrow\tau_k>\tau_{n-1}\\[5pt]
&\Rightarrow\tau_k\in\rrbracket
\tau_{n-1},\infty\rrbracket.
\end{align*}
Conversely, assume that $\tau_k\in\rrbracket \tau_{n-1},\infty\rrbracket$,
i.e., $\tau_k>\tau_{n-1}=\sigma_{n-1}-\sigma_n$. Suppose that
$\sigma_k\notin \rrbracket \sigma_{n-1},\sigma_{n}\rrbracket$,
namely, $\sigma_n<\sigma_k<\sigma_{n-1}$. Then we have
\[\tau_k=\sigma_k-\sigma_n<\sigma_{n-1}-\sigma_n=\tau_{n-1},\]
a contradiction. This yields $\sigma_k\in \rrbracket
\sigma_{n-1},\sigma_{n}\rrbracket$. Similarly, one can verify the case $\sigma_{n-1}<\sigma_n$. This completes the proof. \qed

The following corollary  provides an alternative way to compute the
cyclic  major code.

\begin{coro}\label{tuilun1}
For any permutation $\sigma=\sigma_1\sigma_2\cdots\sigma_n\in S_n$,
define
$$\C(\sigma)=C^{\sigma_n}(\sigma_1\sigma_2\cdots\sigma_{n-1})$$ and
define $L(\sigma)=\sigma_n$. Then we have
\[s_i(\sigma)=i-L(\C^{n-i}(\sigma)),\]
for $1\leq i \leq n$, where $\C^0(\sigma)=\sigma$ and $\C^k(\sigma)=\C(\C^{k-1}(\sigma))$.
\end{coro}
\pf First we see that $\C^{n-i}(\sigma)\in S_i$. By the definition of
$s_n(\sigma)$, we  find
\[ s_n(\sigma)=\#\{\sigma_n+1,\cdots,n\}
=n-\sigma_n=n-L(\sigma)=n-L(\C^0(\sigma)).\]
By the proof of Theorem 1.1, we deduce that
\[ M(\C^{n-i}(\sigma))=(s_1(\sigma),\cdots,s_i(\sigma)),\] which implies that  $s_i(\sigma)=i-L(\C^{n-i}(\sigma))$ for $i=1,2,\cdots,n$. \qed

 The  sequence
\[L(\C^{n-1}(\sigma)),L(\C^{n-2}(\sigma)),\ldots,L(\C^{0}(\sigma))\]
gives an alternative way to compute the cyclic major code. It also
facilitates
 the computation of $H(\sigma)$.  For example, let $\sigma=392648517$. We have
\[M(\sigma)=(0,0,1,3,1,4,3,5,2), \quad H(\sigma) = 496182537, \]
 see Table \ref{xx}.

\begin{table}
\begin{center}
 \begin{tabular}{|c|c|}
   \hline
 & \\
   $\sigma=392648517$   &$I(H(\sigma))=(0,0,1,3,1,4,3,5,2)$\\[5pt]
   $\downarrow$   & $\Uparrow$\\[5pt]
   $\C^0(\sigma)=39264851\textbf{7}$ & $ C_7^{-1}(48617253)\cdot
7=496182537  $\\[5pt]
$\downarrow$   & $\uparrow$\\[5pt]
   $\C^1(\sigma)=5248617\textbf{3}$ & $C_3^{-1}(3751624)\cdot
3=48617253$\\[5pt]
 $\downarrow$   & $\uparrow$\\[5pt]
   $\C^2(\sigma)=271536\textbf{4}$ & $C_4^{-1}(364152)\cdot 4=3751624$\\[5pt]
 $\downarrow$   & $\uparrow$\\[5pt]
$\C^3(\sigma)=53416\textbf{2}$&$C_2^{-1}(25314)\cdot 2=364152$\\[5pt]
$\downarrow$   & $\uparrow$\\[5pt]
$\C^4(\sigma)=3125\textbf{4}$ & $ C_4^{-1}(2431)\cdot
4=25314$\\[5pt]
$\downarrow$   & $\uparrow$\\[5pt]
$\C^5(\sigma)=423\textbf{1}$&$C_1^{-1}(132)\cdot 1=2431$ \\[5pt]
$\downarrow$   & $\uparrow$\\[5pt]
$\C^6(\sigma)=31\textbf{2}$& $C_2^{-1}(12)\cdot 2=132$\\[5pt]
$\downarrow$   & $\uparrow$\\[5pt]
$\C^7(\sigma)=1\textbf{2}$&$C_2^{-1}(1)\cdot 2=12$\\[5pt]
$\downarrow$   & $\uparrow$\\[5pt]
$\C^8(\sigma)=\textbf{1}$&$C_1^{-1}(\emptyset )\cdot 1=1$\\[5pt]
$\downarrow$   & $\uparrow$  \\[5pt]
 (1,2,2,1,4,2,4,3,7)&$\emptyset$ \\[5pt]
$\Downarrow$   &   \\[5pt]
 $M(\sigma)$=(0,0,1,3,1,4,3,5,2)& the construction of $H(\sigma)$ \\[5pt]
&\\
   \hline
 \end{tabular}
\caption {The procedures to compute $H(\sigma)$ and $M(\sigma)$}\label{xx}
\end{center}
\end{table}
The following corollary shows that Han's bijection $H$ commutes with the complementation operator $c$, a property also satisfied by Foata's partial maps and thus by Foata's map $\Phi$.  For a permutation $\sigma\in S_n$, we define
$c\sigma$ as $\tau_1\tau_2\cdots\tau_n$, where $\tau_i=n+1-\sigma_i$. For a
code $a=(a_1, a_2, \ldots, a_n)\in E_n$, we define $ca=(b_1, b_2, \ldots, b_n)$, where $b_i=i-1-a_i$.

\begin{coro}\label{tuilun2}
For $\sigma\in S_n$ and $s=(s_1,s_2,\cdots,s_n)\in E_n$, we have
\[H(c\sigma)=cH(\sigma).\]
\end{coro}

The above corollary can be easily verified by induction on $n$. It
also follows from Theorem \ref{th1} and the relations
\begin{align*}
M(c\sigma)&=c(M(\sigma)),\\[5pt]
I(c\sigma)&=c(I(\sigma)).
\end{align*}

\section {A characterization of fixed points }

In this section, we give a characterization of the fixed points of Han's map
$H$. As will be seen, the fixed points of Han's map are related to the
strong fixed points of Foata's second  fundamental transformation which are easier to
characterize.

 The notion of strong fixed points of Foata's map is related to the strong Foata classes introduced by Bj\"{o}rner and Wachs \cite{bjoner}. A  labeling $w$ of a poset $P$ is called recursive  if every
principal order ideal of $P$ is labeled by a set  of consecutive
 numbers. In particular, if $P$ is a chain
with $n$ elements and $w:P\longrightarrow[n]$ is a labeling of $P$.
Reading the labels from bottom to top, the labels form a permutation
$\sigma=\sigma_1\sigma_2\cdots\sigma_n\in S_n$.
It is easily seen that
$w$ is a recursive labeling of $P$ if and only if  for each $i\in [n]$, $\{\sigma_1,\sigma_2,\cdots,\sigma_i\}$ forms a set of
consecutive numbers.  By the Theorem 4.2 in \cite{bjoner},  we  deduce that a permutation $\sigma\in S_n$  is a strong
fixed point of Foata's map if and only if  for each $i\in [n]$, $\{\sigma_1,\sigma_2,\cdots,\sigma_i\}$ forms a set of consecutive
numbers. For example, $\sigma=45367281\in S_8$ is a
strong  fixed point of Foata's map, while $\pi=34125678$ is not, since
$\{\pi_1,\pi_2,\pi_3\}=\{1,3,4\}$ is not a set of consecutive
numbers.

\begin{theo}\label{th2}
For each  $\sigma\in S_n$, $\sigma$ is  a fixed point of $H$, i.e.
$H(\sigma)=\sigma$, if and only if $\sigma$ is a strong fixed point of Foata's map.
\end{theo}

\pf Suppose that $H(\sigma)=\sigma$. By Theorem \ref{th1}, we see that $I(\sigma)=M(\sigma)$. In particular, we have $s_{n-1}(\sigma)=t_{n-1}(\sigma)$. If
$\sigma_{n-1}>\sigma_n$, by Corollary \ref{tuilun1} we have
\[ s_{n-1}(\sigma)=n-1-L(\C(\sigma))=
n-1-(\sigma_{n-1}-\sigma_n)=n-1+\sigma_n-\sigma_{n-1},\]
and  by definition $t_{n-1}(\sigma)=n-\sigma_{n-1}$. It follows that
$\sigma_n=1$. If $\sigma_{n-1}<\sigma_n$, then
$s_{n-1}(\sigma)=\sigma_n-\sigma_{n-1}-1$ and
$t_{n-1}(\sigma)=n-\sigma_{n-1}-1$. Hence  $\sigma_n=n$. Using   relation \eqref{eq1}, we get
\begin{equation}\label{mc} M(\C(\sigma))=(s_1(\sigma),s_2(\sigma),\cdots,s_{n-1}(\sigma)).
\end{equation}
Moreover, when $\sigma_n=1$ or $\sigma_n=n,$ we have
\begin{equation}\label{cs}
\C(\sigma)=C_{\sigma_n}(\sigma_1\cdots\sigma_{n-1}).
\end{equation} Combining (\ref{mc}), (\ref{cs}) and the fact that \[I(C_{\sigma_n}(\sigma_1\cdots\sigma_{n-1}))
=(t_1(\sigma),t_2(\sigma),\cdots,
t_{n-1}(\sigma)),\] we obtain
\[M(\C(\sigma))=I(\C(\sigma)).\]
   By induction, we deduce that $\C(\sigma)$ is a strong fixed point of Foata's map. Consequently, by  relation (\ref{cs}), we have
\[\{\sigma_1,\sigma_2,\cdots,\sigma_i\}
=\left\{
                                        \begin{array}{ll}
                                          \{(\C(\sigma))_1+1,(\C(\sigma))_2+1,\cdots,(\C(\sigma))_i+1\}, & \hbox{if $\sigma_n=1;$} \\[5pt]
                                          \{(\C(\sigma))_1,(\C(\sigma))_2,\cdots,(\C(\sigma))_i\}, & \hbox{if $\sigma_n=n,$}
                                        \end{array}
                                      \right.\]
which is a set of consecutive integers. Thus $\sigma$ is a strong  fixed point of Foata's map.

Conversely, suppose that $\sigma\in S_n$ is a strong  fixed point of Foata's map.
So $\{\sigma_1,\cdots,\sigma_{n-1}\}$ is a set of consecutive
integers with $n-1$ numbers in $[n]$. This implies that $\sigma_n=1 $
or $\sigma_n=n$. Hence
\[C^{\sigma_n}(\sigma')=C_{\sigma_n}(\sigma')=\left\{
               \begin{array}{ll}
                 (\sigma_1-1)\cdots(\sigma_{n-1}-1), & \hbox{if $\sigma_n=1$;} \\[5pt]
                \sigma_1\cdots\sigma_{n-1}, & \hbox{if $\sigma_n=n$,}
               \end{array}
             \right.
\]
where $\sigma'=\sigma_1\sigma_2\cdots\sigma_{n-1}$. It follows that $C^{\sigma_n}(\sigma')$ is a
strong fixed point of Foata's map. By induction, we deduce that \[ H(C^{\sigma_n}(\sigma'))=C^{\sigma_n}(\sigma').\] Hence
\begin{align*}
H(\sigma)&=C_{\sigma_n}^{-1}(H(C^{\sigma_n}(\sigma')))\cdot\sigma_n\\[5pt]
&=C_{\sigma_n}^{-1}(C^{\sigma_n}(\sigma'))\cdot\sigma_n \\[5pt] &=C_{\sigma_n}^{-1}(C_{\sigma_n}(\sigma'))\cdot\sigma_n=
\sigma'\cdot\sigma_n=\sigma,
\end{align*}
as desired. This completes the proof.  \qed

The following corollary gives another characterization of the  fixed points of $H$ in terms of the cyclic major code and the inversion code.

\begin{coro} \label{tuilun3}
Let $\sigma\in S_n$. The following statements are equivalent:
\begin{itemize}
\item[(1)] $M(\sigma)=I(\sigma)$, that is, $\sigma$ is a fixed point of $H$. \\
\item[(2)] $I(\sigma)=(t_1(\sigma),t_2(\sigma),\cdots,t_n(\sigma))$ such that $t_i(\sigma)=0$ or $i-1$ for each $i\in [n]$.
\end{itemize}
\end{coro}

\pf It is easy to check that $\sigma$ satisfies Condition (2) if $\sigma$ is a strong  fixed point of Foata's map.
 Conversely, suppose that $I(\sigma)=(t_1(\sigma),t_2(\sigma),\cdots,t_n(\sigma))$ such that $t_i(\sigma)=0$ or $i-1$ for each $i\in [n]$.
We proceed by induction on $n$ to show that $\sigma$ is a strong  fixed point of Foata's map.
The statement is obvious for $n=1$.
Now we  assume that the assertion  holds for any permutation of length $n-1$
satisfying Condition (2). It is clear that \[ I(C_{\sigma_n}(\sigma'))=(t_1(\sigma),\ldots,t_{n-1}(\sigma)).\]
The inductive hypothesis implies that $I(C_{\sigma_n}(\sigma'))$ is a strong  fixed point of length $n-1$.
Since $t_n=0$ or $t_n=n-1$,  we have $\sigma_n=1$ or $\sigma_n=n$, and hence
\[ C_{\sigma_n}(\sigma')=\left\{
               \begin{array}{ll}
                 (\sigma_1-1)\cdots(\sigma_{n-1}-1), & \hbox{if $\sigma_n=1$;} \\[5pt]
                \sigma_1\cdots\sigma_{n-1}, & \hbox{if $\sigma_n=n$.}
               \end{array}
             \right.
\] It follows that $\sigma$ is also a strong  fixed point of Foata's map. Now the corollary is a consequence of  Theorem \ref{th2}.\qed

\begin{coro} \label{tuilun3}
For any $n\geq1$, Han's map $H$  has
 $2^{n-1}$ fixed points.
\end{coro}

By Theorem \ref{th2}, we see that every
fixed point of $H$ is a fixed point of $\Phi$, but the converse
is not true. For example, let $\sigma=14235\in S_5$. Then $\sigma$ is a fixed point of $\Phi$, but it is not a fixed point of $H$.

\vspace{0.5cm}
 \noindent{\bf Acknowledgments.} We wish to thank the referee for valuable suggestions. This work was supported by  the 973
Project, the PCSIRT Project of the Ministry of Education,  and the National Science
Foundation of China.

\end{document}